\documentstyle[twoside]{article}
\input amssym.def
\input scrload.tex
\catcode`@=11
\long\def\@makefntext#1{\noindent #1}
\newskip\tabcentering \tabcentering=1000pt plus 1000pt minus 1000pt
\def\REF#1{\par\hangindent\parindent\indent\llap{#1\enspace}\ignorespaces} 
\def\MCH#1#2{\setbox0=\hbox{\raise#1\hbox{#2}}\smash{\box0}}

\def\@evenfoot{}\def\@oddfoot{}

\newtheorem{thm}{Theorem}[section]
\newtheorem{cor}[thm]{Corollary}

\newtheorem{lem}[thm]{Lemma}

\newtheorem{defn}[thm]{Definiteion}
\newtheorem{rem}[thm]{Remark}

 \newcommand{\h}{\mathcal{H}}
 
\newcommand{\mC}{\mathcal{C}}

\newcommand{\mP}{\mathcal{P}}

 \newcommand{\M}{\mathcal{M}}

  \newcommand{\mZ}{\mathcal{Z}}
 \newcommand{\Y}{\mathcal{Y}}
 
\newcommand{\N}{\mathbb{N}}

\catcode`@=12

\def\N{{\Bbb N}}

\def\C{{\Bbb C}}

\def\bc{\begin{center}}
\def\ec{\end{center}}
\def\no{\noindent}
\def\hang{\hangindent\parindent}
\def\textindent#1{\indent\llap{\qquad #1\ \ \enspace}\ignorespaces}
\def\ref{\par\hang\textindent}

plus 1pt minus 1pt
\vspace{2 true cm}
\begin{document}
\bc{\large\bf Relations between positive definite functions and irreducible representations on a locally compact groupoid
\\
\footnotetext{\footnotesize}}
 \ec

\vspace*{0.1 true cm}
\bc{\bf H. Amiri\\
{\small\it Department of Mathematics, University of Isfahan, Isfahan, Iran\\
E-mail: h.amiri@sci.ui.ac.ir}}\ec
 \vspace*{3 true mm} \noindent{\small
{\small\bf Abstract}\\
If $G$ is a locally compact groupoid with a Haar system $\lambda$, then a positive definite function $p$ on $G$ has a form $p(x)=\langle L(x)\xi(d(x)),\xi(r(x))\rangle$, where $L$ is a representation of $G$ on a Hilbert bundle ${\h}=(G^0,\{H_u\},\mu)$, $\mu$ is a quasi invariant measure on $G^0$ and $\xi\in L^{\infty}(G^0,{\h})$. [10].

In this paper firt we prove that if $\mu$ is a quasi
invariant ergodic measure on $G^0$, then two
corresponding representations of $G$ and $C_c(G)$ are irreducible in the same time. Then by using the theory of positive linear functionals on $C^*(G)$ we show that when $\mu$ is an ergodic quasi invariant measure on $G^0$, for a positive definite function $p$ which is an extreme point of ${\mP}^{\mu}_1(G)$ the corresponding representation $L$ is irreducible and conversely, every irreducible representation $L$ of $G$ on a Hilbert bundle ${\h}=(G^0,\{H_u\},\mu)$ and every section $\xi\in {\h}(\mu)$ with norm one, define an extreme point of ${\mP}^{\mu}_1(G)$.

\vspace{3mm}

\no{\small\bf Keywords:} Topological groupoid, Irreducible
representation, positive definite functions. \ \

\no{\small\bf MR(2000) Subject Classification}\ \ {\rm 22A22,
22A25, 43A35}}
\vspace{.5cm}
\baselineskip 15pt\\
{\bf Introduction}\\
In the representation theory of locally compact groups, unitary
representations of such a group $G$ are integrated up with respect
to the left Haar measure against $C_c(G)$ (or even $L^1(G))$
functions to give a representation of the convolution algebra
$C_c(G)$. Also all the representations of $L^1(G)$ are obtained in this
way and two corresponding representations of
$G$ and $C_c(G)$ are irreducible at the same time. [4, section.22].

In the locally compact groupoid case, a representation of such a groupoid is defined by a Hilbert bundle. Jean Renault in
 [8] showed that, as in the group case, integrating up such a representation with respect to the left Haar system gives a
 $*-$representation of $C_c(G)$, the Hilbert space of the representation of $C_c(G)$ is the space of square integrable
  section of the Hilbert bundle. Renault also showed that all representations of $C_c(G)$ are obtained in this way.

The set of positive definite functions on a locally compact groupoid $G$ with a left Haar system $\lambda$ is defined by
Arlan Ramsay and Martin E. Walter in [10] and for measured groupoids is by Jean Renault in [9]. In this paper the set of positive definite function of a measured groupoid $(G,\lambda,\mu)$ is denoted by ${\mP}^{\mu}(G)$.
If $(G,\lambda,\mu)$ is a measured groupoid and $L$ is a representation of $G$ on a Hilbert bundle ${\h}$ and  $\xi$ is a section in $L^{\infty}(G^0,{\h})$, then $p(x)=\langle L(x)\xi(d(x)),\xi(r(x))\rangle$ defines an element of ${\mP}^{\mu}(G)$.  Conversely if $p\in {\mP}^{\mu}(G)$, then there exists a representation $L$ of $G$ on a Hilbert bundle ${\h}$ and a section $\xi\in L^{\infty}(G^0,{\h})$ such that $p(x)=\langle L(x)\xi(d(x)),\xi(r(x))\rangle$ [9, Theorem 1.1] and [10].

In Section $2$ of this paper we prove that if $\mu$ is a quasi
invariant ergodic measure on $G^0$, then two
corresponding representations of $G$ and $C_c(G)$ are irreducible in the same time.
In section $3$ for quasi invariant measure $\mu$ we introduce a subset of ${\mP}^{\mu}(G)$ which is denoted by
${\mP}_1^{\mu}(G)$. After proving some necessary lemma and by using the theory of positive
linear functionals on $C^*(G)$, when $\mu$ is an ergodic quasi invarianr measure we
prove that if $p$ is an element of ${\mP}_1^{\mu}(G)$ and is
an extreme point of ${\mP}_1^{\mu}(G)$, then the corresponding
representation to $p$ is irreducible, and conversely if $\mu$ is an ergodic quasi invariant measure on $G^0$ and $L$
is an irreducible representation of $G$ on a Hilbert bundle
${\h}=(G^0,\{H_u\},\mu)$ and $\xi\in{\h}(\mu)$ is a square integrable section with
norm one, then $p(x)=\langle L(x)\xi(d(x)), \xi(r(x))\rangle$ for
$x\in G$ is an extreme point of ${\mP}_1^{\mu}(G)$.

\section{Preliminaries}
The reader who is unfamiliar with {\it groupoids} will find the definiteion of a groupoid and also topological, locally compact groupoid in [6], [8]. Throughout the paper we shall use the definiteion of a locally compact
groupoid given by Alan L.T.Paterson in [6, P.31].

We just indicate that the range and domain
maps on $G$ are $r:G\rightarrow G^0$, $d: G \rightarrow G^0$ by setting
$r(x)=xx^{-1},\ d(x)=x^{-1}x$ and $G^0=r(G)=d(G)$ is the unit
space of $G$. Each of the maps $r,d$ fibers the groupoid $G$
over $G^0$ with fibers $\{G^u\}$, $\{G_u\}\ \ (u\in G^0)$, so that
$G^u=r^{-1}(\{u\}),\ G_u=d^{-1}(\{u\})$.

A left Haar system for a locally compact groupoid $G$ is a family
$\lambda=\{\lambda^u\}\ (u\in G^0)$, where each $\lambda^u$ is a positive
regular Borel measure on $G^u$, with support equal to $G^u$. In
addition, for any $g\in C_c(G)$, the function $g^0$, where
$g^0(u)=\int_{G^u}g(z) d\lambda^u(z)$ belongs to $C_c(G^0)$ and
for any $x\in G$ and $f\in C_c(G)$,
$$\int_{G^{r(x)}}f(x^{-1}z)d\lambda^{r(x)}(z)=\int_{G^{d(x)}}f(y)d\lambda^{d(x)}(y).\ \ \ (*)$$

Two positive measures $\mu$ and
$\mu'$ are called {\it equivalent} if each one is absolutely
continuous with respect to the other, and in this case we write
$\mu\sim\mu'$.

For any topological space $X$ we denote by $B(X)$ the $\sigma-$
algebra of Borel subsets of $X$.

For a locally compact groupoid $G$ with a left Haar system $\lambda$ if $\mu$ is a probability measure on $G^0$ then $\mu$ and
$\lambda$ determine two regular Borel positive measures
$\nu$ and $\nu^{-1}$ on $B(G)$  with
$\nu=\int_{G^0}\int_{G^u}\lambda^ud\mu(u)$ and
$\nu^{-1}(E)=\nu(E^{-1}) (E\in B(G))$, and if $\mu_1\sim\mu_2$ then $\nu_1\sim\nu_2$.
Also the measure $\nu^2$ on $G^2$ is defined by
$\nu^2=\int_{G^0}\int_{G^u}\int_{G_u}d\lambda_u(x)d\lambda^u(y)d\mu(u)$, where $\lambda_u=(\lambda^u)^{-1}$ on $G_u$.
The probability measure $\mu$ is called {\it quasi invariant} if
$\nu\sim\nu^{-1}$. If $\mu$ is quasi invariant, then the
Radon-Nikodym derivative $\delta_{\mu}=\frac{d\nu}{d\nu^{-1}}$ is called the
modular function of $\mu$ and $d\nu_0=\delta_{\mu}^{\frac{-1}{2}}d\nu$
defines a measure which satisfies $\nu_0(E)=\nu_0(E^{-1})$ for
every Borel subset $E$ of $G$.( cf. [6, p. 87]).

We denoted by $Q$ the set of all quasi invariant measure on $G^0$. If $\mu\in Q$ and $U\subseteq G^0$ is a $\mu-$conull Borel set, then $G|_U=r^{-1}(U)\cap d^{-1}(U)$, which is a subgroupoid and $\nu-$conull, is called {\it inessential reduction} of $G$.

Throughout the paper we shall  use the theory of Hilbert bundles
given by Dixmier in [2, Part II].

\begin{defn}  {\rm A  {\it Hilbert bundle} ${\h}$ is a triple such as $(X,
\{H_u\},\mu)$. Here $X$ is a (second countable) locally compact
Hausdorff space and $\mu$ is a probability measure on $X$ and
$H_u$ is a Hilbert space for every $u\in X$. A {\it section} of a
Hilbert bundle ${\h}$ is a function $f:X\rightarrow
\bigcup_{u\in X}H_u$ where $f(u)\in H_u$. There is also a linear
subspace $S$ consisting of sections which are called $\mu-$measurable
 and there exists a sequence $(\xi_1, \xi_2, \xi_3,\ldots )$ of
elements of $S$ such that, for every $u\in X$,
$\{\xi_n(u)\}_{n=1}^{\infty}$ forms a total sequence in $H_u$.
The sequence $\{\xi_i\}_{i=1}^{\infty}$ is called a {\it
fundamental sequence} for ${\h}$}.
\end{defn}

The Hilbert space
${\h}(\mu)=L^2(X,\{H_u\},\mu)$ is defined in the obvious way as the space
of equivalence classes of measurable sections $f$ for which the
function $u\rightarrow \|f(u)\|_2^2$ is $\mu-$ integrable. The space ${\h}(\mu)$ with the inner product
$$\langle f,g\rangle =\int_X\langle f(u),g(u)\rangle d\mu(u)\quad (f, g\in{\h}(\mu)$$ is a Hilbert space [2, P.168].

For two Hilbert spaces  $H$ and $K$, an  operator $T\in B(H,
K)$ ( the space of all bounded linear operators from $H$ into
$K$) is called {\it unitary} if $T$ is bijective and isometry.
The space of all unitary operators in $B(H, K)$ is denoted by
$U(H, K)$.

\begin{defn} {\rm [5, P. 69]  An a.e. representation of the locally compact groupoid
$G$ consists of a $\mu\in Q$ (with the associated measures $\nu$ and $\nu^2$) and a Hilbert bundle ${\h}=(G^0 ,\{H_u\} ,\mu)$ and a $\mu-$conull subset $U$ of $G^0$ with $L(x)\in U((H_{d(x)}, H_{r(x)})$ for  $\nu.\ a.e.\ x\in G|_U$ such that:\\
1) $L(u)=I(u)$ is the identity operator on $H_u$ for $\mu.\ a.e\ u\in G^0$,\\
2) $L(x)L(y)=L(xy)\ \ \ \nu^2.\ a.e\ (x,y)\in G^2$, \\
3) $L(x)^{-1}=L(x^{-1})\ \ \  \nu.\ a.e\ x\in G $,\\
4) for any $\xi,\eta\in{\h}(\mu)$, the function
$x\rightarrow\langle L(x)\xi(d(x)),\eta(r(x))\rangle$ is $\nu-$
measurable on $G$}.
\end{defn}
The space of continuous functions with compact support is denoted by
$C_c(G)$. The space $C_c(G)$ with $\|.\|_I$ and involution
$f^*(x)=\overline{f(x^{-1})}$ is a normed $*-$ algebra [6, P.40]. If $L$
is a representation of the locally compact groupoid $G$ on a Hilbert
space ${\h}$, then $\pi_L: C_c(G)\rightarrow
B({\h}(\mu))$ given by :
$$\langle \pi_L(f)\xi,\eta\rangle=\int_Gf(x)\langle L(x)(\xi(d(x)),\eta(r(x))\rangle d\nu_0(x)$$
is a representation of $C_c(G)$ of norm $\leq 1$. Also every
representation of $C_c(G)$ is of the form $\pi_L$ for some
representation $L$ of $G$ and the correspondence $L\rightarrow\pi_L$
preserves the natural equivalence relation on the representation of
$G$ and the representation of $C_c(G)$ [6, P. 98] and [8, P. 72]. The completion of $C_c(G)$ with respect to
$$\|f\|=\sup_{L}\|\pi_L(f)\|$$
is defined to be ([8, P.58]) the $C^*-$algebra of $G$ denoted $C^*(G)$.

Let $\mu\in Q$ and  ${\h}=(G^0,\{H_u\},\mu)$ be a Hilbert bundle, and $T$ be a section $u\rightarrow T(u)\in B(H_u)$ of bounded linear operators, then the section $T$ is called measurable if for each measurable section $\xi$, the section $T\xi$ defined by $(T\xi)(u)=T(u)\xi(u)$ is measurable. Also the section $T$ is called essentially bounded if the essential sup $u\rightarrow \|T(u)\|$ is finit. A linear operator $T: {\h}(\mu)\rightarrow {\h}(\mu)$ is called decomposable if it is defined by an essentially bounded measurable section of bounded operators $u\rightarrow T(u)\in B(H_u)$. In this case we write $T=\int^\oplus T(u)d\mu(u)$, where
$$\langle T\xi,\eta\rangle=\int_{G^0}\langle T(u)(\xi(u)),\eta(u)\rangle d\mu(u)\ \ \ (\xi,\eta\in{\h}(\mu).$$

For $\mu\in Q$ let $L_{\C}^{\infty}(G^0,\mu)$ be the set of all Complex valued $\mu-$measurable essential bounded function on $G^0$, in which we identify two functions that are equal almost everywhere, and $C_0(G^0)$ be the set of complex-valued function on $G^0$ which are continuous and vanish at infinity. Also $(G^0,\{H_u\},\mu)$ a Hilbert bundle, then if $f\in L_{\C}^{\infty}(G^0,\mu)$ or if $f\in C_0(G^0)$, the section of operators $u\rightarrow f(u)I(u)\in B(H_u)$ is a measurable and essentially bounded, where $I(u)$ is the identity operator of $H_u$. We denote by $T_f$ the corespondence operator on $L^2(G^0,\{H_u\},\mu)$ which is $T_f=\int^{\oplus}f(u)I(u)d\mu(u)$.

\begin{defn} an operator $T\in B(L^2(G^0,\{H_u\},\mu))$ is called diagonalisable [resp. continuously diagonalisable] if $T$ has the form $T_f$ for some $f\in L_{\C}^{\infty}(G^0,\mu)$ [resp. $f\in C_0(G^0) ]$. We put
$${\mZ}=\{T_f: f\in L_{\C}^{\infty}(G^0,\mu)\}$$
$${\Y} =\{T_f: f\in C_0(G^0)\}$$
\end{defn}
It is easy to check that if $T\in B(L^2(G^0,\{H_u\},\mu))$ is a decomposable operator, then $TU=UT$ for each $U\in {\mZ}$.

\begin{lem} [Dixmear page 187] Let $\mu\in Q$ and  $(G^0,\{H_u\},\mu)$ be a  Hilbert bundle and $T:L^2(G^0,\{H_u\},\mu)\rightarrow  L^2(G^0,\{H_u\},\mu)$ be a bounded linear operator. If $TU=UT$ for each $U\in{\mZ}$, then $T$ is a decomposable operator.
\end{lem}
\begin{cor} [Dixmier, P.188]
For an operator $T\in B(L^2(G^0,\{H_u\},\mu))$ to be decomposable it is necessary and sufficient that it commute with the diagonalisable operators.
\end{cor}
Since for a locally compact groupoid $G$, the set of units, $G^0$ is $\sigma-$ compact set therefor by [Dixmier, P.188] we have the following Lemma.
\begin{lem}
If $T:L^2(G^0,\{H_u\},\mu)\rightarrow  L^2(G^0,\{H_u\},\mu)$ is a bounded linear operator and $TU=UT$ for each $U\in {\Y}$, then $T$ is a decomposable operator.
\end{lem}

\section{Hilbert bundle with ergodic measure on $G^0$}
We start this section by the following definiteion.

\begin{defn} {\rm Let $\mu\in Q$. A measurable set $A$ in $G^0$ is {\it almost invariant} (with respect to $\mu$) if for $\nu.\ a.e\  x\in G$, $r(x)\in A$ if and only if $d(x)\in A$. The measure $\mu$ is called {\it ergodic} if every almost invariant measurable set is null or conull [8, P.24].}
\end{defn}
We refer the reader for the definiteion of Hilbert subbundle to [2, P 173].
\begin{defn} {\rm Let $\mu\in Q$ and $L$ be a representation of $G$ on a Hilbert bundle ${\h}=(G^0,\{H_u\},\mu)$
, then a Hilbert subbundle ${\M}=(G^0,\{M_u\},\mu)$ of ${\h}$ is {\it nontrivial invariant} if:

1) $L(x)M_{d(x)}\subseteq M_{r(x)}\ \ \ \nu.\ a.e\ x\in G$,

2) $\mu\{u: M_u\neq 0\}>0$ and $\mu\{u: M_u\neq H_u\}>0$.}
\end{defn}
\begin{rem} {\rm When ${\M}$ is a Hilbert subbunle, then the section $u\rightarrow P(u)\in B(H_u)$
 is a bounded  measurable section of projections, where $P(u)$ is the projection corresponding to $M_u$ [2, P. 173].
 Working with fundamental sequence we can show that if ${\M}$ is a nontrivial Hilbert subbundle,
  then there exists two nonzero sections $\xi, \eta\in{\h}(\mu)$, which $\xi(u)\in M_u, \eta(u)\in M_u^{\perp}\ \ (\mu\  a.e\ u\in G^0)$.
   By [2, P. 173], ${\M}(\mu)=L^2(G^0,\{M_u\},\mu)$ is a closed subspace of ${\h}(\mu)$ and by [2, P. 183] $P=\int^{\oplus}P(u)d\mu(u)$
    is the projection corresponding to ${\M}(\mu)$.}
\end{rem}

\begin{defn} {\rm For a locally compact groupoid $G$ with a left Harr system $\lambda$, a representation $L$ on a Hilbert bundle ${\h}$
 is called {\it reducible} if $L$ admits a nontrivial invariant Hilbert subbundle. And $L$ is called {\it irreducible} if it is not reducible.
 A representation $\pi$ of $C^*(G)$ on a Hilbert space $H$ is called irreducible if there isn't any nontrivial invariant closed subspace of $H$ under
 $\pi$.}
\end{defn}

\begin{rem}
{\rm If $\mu$ is a non ergodic quasi invariant measure, then there
exists an almost invariant measurable subset $E\subseteq G^0$ with
$\mu(E)> 0$ and $\mu(E^c)> 0$. Now if $L$ is a representation of $G$
on a Hilbert bundle $(G^0,\{H_u\},\mu)$, then
$(G^0,\{M_u\},\mu)$ is a Hilbert subbundle, where
$$M_u=
\left\lbrace
  \begin{array}{ll}
H_u&\mbox{if}\ u\in E\\
0&\mbox{otherwise,}
  \end{array}\right.$$
and is obviously nontrivial. It is easy to check that it is invariant under $L$.
 Therefore when we talk about an irreducible representation, $\mu$ must be an ergodic measure. Also when $\mu$ is an ergodic measure, then a Hilbert
  subbundle $(G^0,\{M_u\},\mu)$ of $(G^0,\{H_u\},\mu)$ is  nontrivial if and only if
$\mu\{u: 0\neq M_u\neq H_u\}>0$}
\end{rem}

We will to show that if $\mu\in Q$ is an ergodic measure and $L$ is a representation of $G$ on a Hilbert bundle ${\h}$, then
$(L, {\h})$ is irreducible iff $(\pi_L, {\h}(\mu))$ is irreducible.

\begin{lem} Let ${\h}=(G^0,\{H_u\},\mu)$ be a Hilbert bundle,  $P\in{\h}(\mu)$ be a decomposable projection and $P=\int T(u)d\mu(u)$,
 then $T(u)$ is a projection for $\mu.\ a.e\ u\in G^0$.
\end{lem}
\begin{proof} Suppose that $M$ is the closed subspace corresponding to $P$. Therefore $\xi\in M$ if and only if $P\xi=\xi$.
 Since $(P\xi)(u)=T(u)\xi(u)$ for $\mu.\ a.e\ u\in G^0$, so $\xi\in M$ if and only if $T(u)\xi(u)=\xi(u)$ for $\mu.\ a.e\ u\in G^0$. Put
 $M_u=$ closed linear span $\{\xi(u): \xi\in M\}$
 and let $\{\xi_n\}_{n=1}^{\infty}$ be a fundamental sequence for ${\h}$. Since $P\xi_n\in M$, so $T(u)(\xi_n(u))\in M_u\ \mu .a.e $. We can find a $\mu-$conull set $F\subset G^0$ with $T(u)(\xi_i(u))\in M_u$ for $u\in F$ and $i\in\N$.  Since $\{\xi_n(u)\}_{n=1}^{\infty}$
  spans a dence subspace of $H_u$, we have $T(u)H_u\subseteq M_u$ for $u\in F$. Now let $u\in F$ and $x\in M_u$, then there exists a sequence $\xi_k(u)$
such that $\xi_k\in M$ and $\xi_k(u)\rightarrow x$.
If $E\subseteq G^0$ be a $\mu-$ conull set which $T(u)\xi_k(u)=\xi_k(u)$ for $u\in E$ and $k=1, 2, 3\ldots$,
then for $u\in E\cap F$, $T(u)\xi_k(u)=\xi_k(u)\rightarrow T(u)x$. So for $u\in E\cap F$, $T(u)x=x$ and therefore $T(u)$ is the projection corresponding to $M_u$ for $\mu.\ a.e\ u\in G^0.$
\end{proof}

\begin{rem} {\rm [8, p.59] If $g$ is a bounded continuous function on $G^0$ and $f\in C_c(G)$, we define $(gf)(x)=g(r(x))f(x)$, then $gf\in C_c(G)$. It is easy to show that if $L$ is a representation of $G$, then $\pi_L(gf)=T_g\pi_L(f)$.}
\end{rem}

\begin{thm} Let $\mu\in Q$ be an ergodic measure and $L$ be a  representation of $G$ on a Hilbert bundle ${\h}=(G^0,\{H_u\},\mu)$. If
$(\pi_L, {\h}(\mu))$ is a reducible representation of $C^*(G)$, then $L$ is a reducible representation.
\end{thm}
\begin{proof} Suppose that $0\neq M\neq {\h}(\mu)$ is a closed subspace of ${\h}(\mu)$ invariant under $\pi_L$. Then $P_M\in (\pi_L(C^*(G)))'$
(the commutant of $\pi_L(C^*(G))$), where $P_M$ is the projection corresponding to $M$. Let $g\in C_0(G^0)$ and
$\{f_{\alpha}\}$ be an approximate identity for $C^*(G)$, therefore
$f_{\alpha}*h\rightarrow h\ \ (h\in C_c(G))$ which implies that
$\pi_L(f_{\alpha}*h)\rightarrow \pi_L(h)$. So
$T_g\pi_L(f_{\alpha})\pi_L(h)\xi=T_g\pi_L(f_{\alpha}*h)\xi\rightarrow
T_g\pi_L(h)\xi$ for each $\xi\in{\h}(\mu)$. Hence
$T_g\pi_L(f_{\alpha})\rightarrow T_g$ in strong operator topology,
by Remark 2.7,
$T_g\in\overline{(\pi_L(C_c(G)))}^{S.O.T}=(\pi_L(C_c(G))^{''}$. So
$P_MT_g=T_gP_M$. Hence by Lemma 1.3 $P_M$ is decomposable, and by Lemma 2.6,
 $P=\int^{\oplus}P(u)d\mu (u)$, where $P(u)$ is a projection in
$B(H_u)$ for all elements $u$ in a $\mu-$ conull set $E$. If we put
$M_u=P(u)H_u$ for $u\in E$, then it is easy to show that
$M=\{\xi\in{\h}(\mu) :\ \xi(u)\in M_u\}$ and
$M^{\perp}=\{\xi\in{\h}(\mu) :\ \xi(u)\in M_u^{\perp}\}$. Since for
every measurable section  $\xi\in{\h}(\mu)$ the section
$u\rightarrow P(u)\xi(u)$ is measurable,
 ${\M}=(G^0,\{M_u\},\mu)$ is a Hilbert bundle of ${\h}$ [2, P. 173]. Also nontriviality of $M$, implies that
${\M}$ is nontrivial. The only thing which must be
proved is $L(x)M_{d(x)}\subseteq M_{r(x)}\ \ \nu.\ a.e$. Let
$\{\xi_n\}_{n=1}^{\infty}$ be a fundamental sequence for Hilbert
bundle ${\h}$, then $\{P\xi_n\}_{n=1}^{\infty}$ and
$\{(I-P)\xi_n\}_{n=1}^{\infty}$ are fundamental sequence for
${\M}$ and ${\M}^{\perp}=(G^0,\{M_u^{\perp}\},\mu)$ respectively.
Now let $f\in C^*(G)$, then for $m, n\in\N$,
$$\begin{array}{ll}
\int_Gf(x)\langle L(x)P(d(x))(\xi_n(d(x))), (I(r(x))-P(r(x)))\xi_m(r(x))\rangle d\nu_0(x)&\\
=\int_Gf(x)\langle L(x)(P\xi_n)(d(x)), ((I-P)(\xi_m))(r(x))\rangle d\nu_0(x)&\\
=\langle \pi_L(f)P\xi_n, (I-P)\xi_m\rangle&\\
=0
\end{array}
.$$
Therefore $\langle L(x)P(d(x))(\xi_n(d(x))), (I(r(x))-P(r(x)))\xi_m(r(x))\rangle=0\ \ \nu.\ a.e$. But $\{(I(r(x))-P(r(x)))\xi_m(r(x))\}_{m=1}^{\infty}$ for $x\in G|_E$ spans a dence subspace of $M_{r(x)}^{\perp}$, therefore $L(x)P(d(x))(\xi_n(d(x)))\subseteq M_{r(x)}$ for $x\in G|_E$. So $L(x)M_{d(x)}\subseteq M_{r(x)}$ for $x\in G|_E$, where $G|_E$ is $\nu-$conull set.
\end{proof}

\begin{thm} Let $\mu\in Q$ be an ergodic measure and $L$ be a reducible representation of $G$ on a Hilbert bundle ${\h}=(G^0,\{H_u\},\mu)$, then
$(\pi_L, {\h}(\mu))$ is a reducible representation of $C^*(G)$.
\end{thm}
\begin{proof} Suppose that ${\M}=(G^0,\{M_u\},\mu)$ is a nontrivial invariant Hilbert subbundle for the representation $L$, then by Remark 2.3
${\M}(\mu)$ is a nontrivial closed subspace of ${\h}(\mu)$.
 We show that ${\M}(\mu)$ is invariant under the representation $\pi_L$. Let $f\in C^*(G),\ \xi\in{\M}(\mu)$ and $\eta\in({\M}(\mu))^{\perp}$.
  Therefore $U_1=\{u: \xi(u)\in M_u\}$ and $U_2=\{u: \eta(u)\in M^{\perp}_u\}$ are two $\mu-$ conull subsets of $G^0$, hence $U=U_1\cap U_2$ is $\mu-$conull subset of $G^0$ which implies $G|_U$ is $\nu$ conull subset of $G$.
   Also $L(x)M_{d(x)}\subseteq M_{r(x)}\ \ \nu.\ a.e $, so
$$\begin{array}{ll}
\langle \pi_L(f)\xi, \eta\rangle&=
\int_Gf(x)\langle L(x)\xi(d(x)), \eta(r(x))\rangle d(\nu_0)(x)\\
&=\int_{G|_U}f(x)\langle L(x)\xi(d(x)), \eta(r(x))\rangle d(\nu_0)(x)\\
&=0.\end{array}$$
That is $\pi_L(f)\xi\in (({\M}(\mu))^{\perp})^{\perp}={\M}(\mu).\ \ \ \ $
\end{proof}
\vspace{.5cm}

\begin{defn} Suppose $L$ and $L'$ are two representations of a locally
compact groupoid $G$ associated with two Hilbert bundles
$(G^0,\{H_u\},\mu)$ and $(G^0,\{H'_u\},\mu)$ respectively, then
we put
$$\begin{array}{ll}
{\mC}(L, L')=\Big\{\ T:u\rightarrow T_u\in B(H_u,H_u')&\mbox{T is essentialy bunded measurable}
\\& \mbox{and}\ L'(x)T_{d(x)}=T_{r(x)}L(x)\quad
\Big\}
\end{array}$$
and ${\mC}(L)={\mC}(L, L).$
\end{defn}
\begin{defn}
Let $\mu\in Q$ and $L$ a representation of $G$ on a Hilbert bundle $(G^0,\{H_u\},\mu)$. Put
$${\mC}(\pi_L)=\{T\in B(L^2(G^0,\{H_u\},\mu)): \pi_L(f)T=T\pi_L(f)\ \ (f\in C_c(G))\}.$$
\end{defn}
\begin{lem}
If $\mu\in Q$ and  $L$ a representation of $G$ on a Hilbert bundle $(G^0,\{H_u\},\mu)$, then
$${\mC}(\pi_L)=\{T: T=\int^{\oplus}T(u)d\mu(u) ,\mbox{where} \ (T(u))\in \C(L)\}.$$
\end{lem}
\begin{proof}
If $(T(u))\in{\mC}(L)$, then by [Dixmear] the operator $T=\int^{\oplus}T(u)d\mu(u)$ is an element of $B(L^2(G^0,\{H_u\},\mu))$. Let $\xi, \eta\in L^2(G^0,\{H_u\},\mu)$ and $f\in C_c(G)$, then
$$\begin{array}{ll}
\langle\pi_L(f)T\xi,\eta\rangle &=\int f(x)\langle L(x)T(d(x))\xi(d(x)),\eta(r(x))\rangle d\mu\circ\lambda(x)\\
&=\int f(x)\langle T(r(x))L(x)\xi(d(x)),\eta(r(x))\rangle d\mu\circ\lambda(x)\\
&=\int f(x)\langle L(x)\xi(d(x)),T(r(x))^*\eta(r(x))\rangle d\mu\circ\lambda(x)\\
&=\langle\pi_L(f)\xi,T^*\eta\rangle\\
&=\langle T\pi_L(f)\xi,\eta\rangle.
\end{array}$$
So $T\in{\mC}(\pi_L).$\\
Conversaly: Let $T\in {\mC}(\pi_L)$, then $T\in(\pi_L(C_c(G)))'$. In the proof theorem [2.8]
 we shown that if $g$ is a bounded continuous function on $G^0$, then $T_g\in (\pi_L(C_c(G)))''$. Therefore $TT_g=T_gT$ for every continuous bounded function on $G^0$, hence by Lemma 1.6 $T$ is decomposable, that is there exists an essential bounded measurable section $u\rightarrow T(u)\in B(H_u)$ whith $T=\int^{\oplus}T(u)d\mu(u)$. It is enough to show that $L(x)T(d(x))=T(r(x))L(x)\ \ \mu\circ\lambda\ x\in G$. Let $f\in C_c(G)$ and $\xi, \eta\in L^2(G^0,\{H_u\},\mu)$, then
$$\begin{array}{ll}
\int f(x)\langle L(x)T(d(x))\xi(d(x)),\eta(r(x))\rangle d\mu\circ\lambda(x) &=\langle\pi_L(f)T\xi,\eta\rangle\\
&=\langle T\pi_L(f)\xi,\eta\rangle\\
&=\langle\pi_L(f)\xi,T^*\eta\rangle\\
=\int f(x)\lbrace L(x)\xi(d(x)),T(r(x))^*\eta(r(x))\rangle d\mu\circ\lambda(x)&\\
=\int f(x)\lbrace T(r(x))L(x)\xi(d(x)),\eta(r(x))\rangle d\mu\circ\lambda(x)&,
\end{array}$$
hence $\langle L(x)T(d(x))-T(r(x))L(x)\xi(d(x)),\eta(r(x))\rangle=0$ for $\mu\circ\lambda\ x\in G$. So
$L(x)T(d(x))=T(r(x))L(x)\ \ \mu\circ\lambda\ x\in G$.
\end{proof}
\section{Positive definite functions and irreducibility}
\begin{defn} {\rm Let $G$ be a locally compact groupoid
with a Haar system $\lambda$ and $\mu\in Q$. We denote by ${\mP}^{\mu}(G)$ the set of all $\nu-$measurable functions $p$ on $G$ which satisfy for every $f\in C_c(G)$,
$$\int\int f(x)\bar{f}(y)p(y^{-1}x)d\lambda^u(x)\lambda^u(y)\geq
0\ \ \ \mu\ a.e\ u\in G^0.$$
We also put $${\mP}(G)=\bigcup_{\mu\in Q}{\mP}^{\mu}(G).$$
An element of ${\mP}(G)$ is called {\it positive definite function}.}
\end{defn}

Two elements $p_1$ and $p_2$ of ${\mP}(G)$ are {\it equivalent} and we write $p_1\sim p_2$, if $p_1\in {\mP}^{\mu_1}(G)$, $p_2\in{\mP}^{\mu_2}(G))$ with $\mu_1\sim\mu_2$ and
$p_1(x)=p_2(x)\rho(d(x))\rho(r(x))\ \ $ $\nu$ a.e $x\in G$, where $\rho$ is a positive function such that $\rho^2(u)=\frac{d\mu_2(u)}{d\mu_1(u)}$ and $\nu=\mu_1 \circ\lambda$. It is easy to check that "$\sim$" is an equivalence relation.

By [9, Thm 1.1] if $\mu\in Q$ and $p\in L^{\infty}(G^0,\mu)$ is an element of ${\mP}^{\mu}(G)$, then there exists a representation $L$ of $G$ on a Hilbert bundle
 $(G^0,\{H_u\},\mu)$ and a section $\xi$ in $L^{\infty}(G^0,{\h})$ such that
$p(x)=\langle L(x)\xi(d(x)), \xi(r(x)) \rangle$. Conversely if
$\mu\in Q$ and $L$ is a representation of $G$ on a Hilbert bundle
$(G^0,\{H_u\},\mu)$ then a section $\xi\in L^{\infty}(G^0,{\h})$ defines an element
of ${\mP}^{\mu}(G)$ by $p(x)=\langle L(x)\xi(d(x)), \xi(r(x))\rangle$.

We are interested to the elements $p$ of ${\mP}(G)$ which are of the form
$p(x)=\langle L(x)\xi(d(x)), \xi(r(x)) \rangle$ for some
representation $L$ of $G$ on a Hilbert boundel ${\h}$
with $\xi\in {\h}(\mu)$.

\begin{defn} {\rm If $\mu\in Q$ and $L$ is a representation of $G$ on a Hilbert bundle ${\h}=(G^0,\{H_u\},\mu)$ we denote
by ${\mP}^{\mu}_1(G)$ the set of all complex valued functions $p$ on $G$ which are of the form $p(x)=\langle L(x)\xi(d(x)), \xi(r(x)) \rangle$, where $\xi\in{\h}(\mu)$ and $\|\xi\|=1$
and we put $${\mP}_1(G)=\bigcup_{\mu\in Q}{\mP}^{\mu}_1(G)$$}
\end{defn}
\begin{lem} Let $\mu\in Q$ and $L$ be a representation of $G$ on a Hilbert bundle ${\h}=(G^0,\{H_u\},\mu)$ and $\xi\in{\h}(\mu)$, then
$p(x)=\langle L(x)\xi(d(x)), \xi(r(x)) \rangle$ defines an element of
${\mP}^{\mu}(G)$, therefore  ${\mP}^{\mu}_1(G)\subseteq {\mP}^{\mu}(G)$.
\end{lem}
\begin{proof}
 {\rm  By definiteion 1.2 $p$ is $\nu-$measurable. Let $\xi\in{\h}(\mu)$ and $f\in C_c(G)$, then by [3]
 $$\pi_L(f)\xi(u)=\int f(x)L(x)\xi(d(x))d\lambda^u(x)\ \ \mu\ a.e\ u\in G^0,$$ where the integral is weakly defined. Therefore
 $$\begin{array}{ll}
\int\int f(x)\bar{f}(y)p(y^{-1}x)d\lambda^u(x)\lambda^u(y)&=\int\int f(x)\bar{f}(y)\langle L(x)\xi(d(x)), L(y)\xi(d(y))\rangle d\lambda^u(x)\lambda^u(y)\\
&=\langle \int f(x) L(x)\xi(d(x))d\lambda^u(x), \int f(y) L(y)\xi(d(y))d\lambda^u(y)\rangle\geq 0 ,
\end{array}$$
for $\mu\ a.e\ u\in G^0$, therefore $p\in{\mP}^{\mu}(G)$.}
\end{proof}

\begin{rem} {\rm Suppose $\mu, L, \xi, p$ are such as Lemma 3.3, we put
$$\phi_p(f)=\langle\pi_L(f)\xi,\xi\rangle=\int f(x)p(x)d\nu_0(x)$$ for $f\in C^*(G)$, then $\phi_p$ is a positive linear functional on $C^*(G)$ and if $p\in{\mP}^{\mu}_1(G)$, then $\phi_p$ is a state on $C^*(G)$.}
\end{rem}

\begin{lem} Let $\phi$ be a state on $C^*(G)$, then there exists a $\mu\in Q$ and a $p\in {\mP}^{\mu}_1(G)$ such that $\phi(f)=\phi_p(f)$ for $f\in C^*(G)$.
\end{lem}
\begin{proof} Since $\phi$ is a positive linear functional on $C^*(G)$, then by GNS construction there exists
 a representation $(\pi_{\phi}, H_{\phi})$ of $C^*(G)$, unique up to unitary equivalence, with a cyclic vector $\xi_{\phi}$ with norm one such that
 $\phi(f)=\langle\pi_{\phi}(f)\xi_{\phi},\xi_{\phi}\rangle$ for $f\in C^*(G)$.
  Since every representation of $C^*(G)$ comes from a representation of $G$, there exists a $\mu\in Q$
   and a representation $L$ of $G$ on a Hilbert bundle ${\h}=(G^0,\{H_u\},\mu)$, such that  $(\pi_{\phi}, H_{\phi})$ is unitary equivalent to
   $(\pi_L, {\h}(\mu))$.
Suppose that $T: {\h}(\mu)\longrightarrow H_{\phi}$ is the unitary operator such that $T\pi_L(f)=\pi_{\phi}(f)T$ for every $f\in C^*(G)$.
If $\xi\in{\h}(\mu)$ is such that $T\xi=\xi_{\phi}$, then

$$\begin{array}{ll}
\phi(f)&=\langle \pi_{\phi}(f)\xi_{\phi}, \xi_{\phi}\rangle \\
&=\langle \pi_{\phi}(f)T\xi, T\xi\rangle\\
&=\langle T^*\pi_{\phi}(f)T\xi, \xi\rangle\\
&=\langle \pi_L(f)\xi, \xi\rangle\\
&=\int_Gf(x)\langle L(x)\xi(d(x)), \xi(r(x))\rangle d\nu_0(x).
\end{array}$$

Now if we put $p(x)=\langle L(x)\xi(d(x)), \xi(r(x))\rangle$, then by Lemma 3.3, $p\in {\mP}^{\mu}(G)$ and $\phi(f)=\phi_p(f)$ for $f\in C^*(G)$. Finally $\|\xi_{\phi}\|=1$ implies that $\|\xi\|=1$, so $p\in {\mP}_1^{\mu}(G)$.
\end{proof}
\begin{lem} Let $p(x)=\langle L(x)\xi(d(x)), \xi(r(x))\rangle$ and $p'(x)=\langle L'(x)\xi'(d(x)), \xi'(r(x))\rangle$ are two elements of ${\mP}(G)$, where $L$ and $L'$ are two representations of $G$ respectively on $(G^0,\{H_u\},\mu)$ and $(G^0,\{H'_u\},\mu')$ and $\xi\in{\h}(\mu)$ and $\xi'\in{\h}(\mu')$, then $\phi_{p}=\phi_{p'}$ if and only if  $p\sim p'$.
So the maping $p\rightarrow \phi_p$ from ${\mP}_1(G)$ onto $\Sigma(C^*(G))$ is one to one on ${\mP}^{\mu}_1(G)$ for each $\mu\in Q$.
\end{lem}
\begin{proof} {\rm Suppose that $\phi_{p}=\phi_{p'}$, therefore for $f\in C^*(G)$
 $$\phi_{p}(f)=\langle \pi_L(f)\xi, \xi\rangle=\langle \pi_{L'}(f)\xi', \xi'\rangle\hspace{1cm}(**),$$ it is well known that the representation associated to a positive linear functional $\phi$ on a $C^*$ algebra is unique
  (up to unitary equivalence), therefore $(\pi_L, {\h}(\mu))$ and $(\pi_{L'}, {\h}(\mu'))$
are unitary equivalent. Hence by [8, Cor 2.1.23] $\mu\sim\mu'$.
Let $\rho$ be a positive function which $(\rho(u))^2=\frac{d\mu'(u)}{d\mu(u)}$ then by [10, P.327] $\rho^2(r(x))=\frac{d\nu'(x)}{d\nu(x)}$ and $\rho^2(d(x))=\frac{d\nu'^{-1}(x)}{d\nu^{-1}(x)}$. Therefore $\rho^2(r(x))\delta_{\mu}(x)=\rho^2(d(x))\delta_{\mu'}(x)$ and consequently $\delta_{\mu}^{\frac{-1}{2}}(x)=\frac{\rho(r(x))}{\rho(d(x))}\delta_{\mu'}^{\frac{-1}{2}}(x)$. By (**) for $f\in C_c(G)$
$$\begin{array}{ll}
\int f(x)p(x)\delta_{\mu}^{\frac{-1}{2}}(x) d\nu(x)&=\int f(x)p'(x) \delta_{\mu'}^{\frac{-1}{2}}(x)d\nu'(x)\\
&=\int f(x)p'(x) \delta_{\mu'}^{\frac{-1}{2}}(x)\rho^2(r(x))d\nu(x)\\
&=\int f(x)p'(x) \frac{\rho(d(x))}{\rho(r(x))}\delta_{\mu}^{\frac{-1}{2}}(x)\rho^2(r(x))d\nu(x)\\
&=\int f(x)p'(x) \rho(d(x))\rho(r(x))\delta_{\mu}^{\frac{-1}{2}}(x)d\nu(x)
\end{array}$$
So $p(x)=p'(x)\rho(d(x))\rho(r(x))\ \ \ \nu\ a.e$, hence $p\sim p'$.}

The converse is strightforward.
\end{proof}
\begin{defn} {\rm Let $\mu,\mu'\in Q$ and $L, L'$ are two unitary equivalent representations of $G$ respectively on two Hilbert bundles $(G^0,\{H_u\},\mu)$ and $(G^0,\{H'_u\},\mu')$. For $(\xi,\xi')\in{\h}(\mu)\times {\h}(\mu')$ we write $\xi\sim\xi'$, if there exists a unitary operator $T\in{\mC}(\pi_L,\pi_{L'})$ with $T\xi=\xi'$, where ${\mC}(\pi_L,\pi_{L'})$ is
the commutant of $\pi_L$ and $\pi_{L'}$}
\end{defn}
\begin{lem}  Let  $L, L'$ are two unitary equivalent representations of $G$ respectively on two Hilbert bundles $(G^0,\{H_u\},\mu)$ and $(G^0,\{H'_u\},\mu')$, where $\mu, \mu'\in Q$.
 If $(\xi,\xi')\in {\h}(\mu)\times {\h}(\mu')$ and $\xi\sim\xi'$, then
$p\sim p'$, where $p(x)=\langle L(x)\xi(d(x)),\xi(r(x))\rangle$ and $p'(x)=\langle L'(x)\xi'(d(x)),\xi'(r(x))\rangle$.
\end{lem}
\begin{proof} {\rm By Lemma 3.3 $p\in {\mP}^{\mu}(G)$ and $p'\in {\mP}^{\mu'}(G)$. Also $\mu\sim\mu'$, since $L, L'$ are  unitary equivalent. Suppose that $T\in {\mC}(\pi_L,\pi_{L'})$ is a unitary operator with $T\xi=\xi'$, then
$$\begin{array}{ll}
\int f(x)p(x)\delta_{\mu}^{\frac{-1}{2}}(x)d\nu(x)&=\langle\pi_L (f)\xi, \xi\rangle\\
&=\langle T^*\pi_{L'} (f)T\xi, \xi\rangle\\
&=\langle \pi_{L'} (f)T\xi, T\xi\rangle\\
&=\langle \pi_{L'} (f)\xi', \xi'\rangle\\
&=\int f(x)p'(x)\delta_{\mu'}^{\frac{-1}{2}}(x)d\nu'(x)\\
&=\int f(x)p'(x) \frac{\rho(d(x))}{\rho(r(x))}\delta_{\mu}^{\frac{-1}{2}}(x)\rho^2(r(x))d\nu(x)\\
&=\int f(x)p'(x) \rho(d(x))\rho(r(x))\delta_{\mu}^{\frac{-1}{2}}(x)d\nu(x)
\end{array}$$
Therefore $p(x)=p'(x)\rho(r(x))\rho(d(x))\ \ \nu\ a.e$, where $\rho^2(u)=\frac{d\mu'(u)}{d\mu(u)}$. So $p\sim p'$.}
\end{proof}
\begin{thm} Let $G$ be a locally compact groupoid with a Haar system $\lambda$ and $\mu\in Q$. If $p\in{\mP}^{\mu}_1(G)$ and $p\in Ext({\mP}^{\mu}_1(G))$ then $\phi_p\in Ext(\Sigma(C^*(G)))$, where $\Sigma(C^*(G))$ is the state space of $C^*(G)$.
\end{thm}
\begin{proof} {\rm By Remark 3.4 $\phi_p\in\Sigma(C^*(G))$. Suppose that there are $\phi_1, \phi_2\in \Sigma(C^*(G)) $ such that $\phi_p=\alpha\phi_1+\beta\phi_2$, where $\alpha,\beta\geq 0$ and $\alpha+\beta=1$. By Lemma 3.5 $\phi_1=\phi_{p_1}$ and $\phi_2=\phi_{p_2}$, where $p_1\in{\mP}^{\mu_1}_1(G)$ and $p_2\in{\mP}^{\mu_2}_1(G)$ and $\mu_1, \mu_2\in Q$.
Suppose that $L_1, L_2$ are two representations of $G$ respectively on two Hilbert bundles $(G^0,\{H_u\},\mu_1)$ and $(G^0,\{H'_u\},\mu_2)$ with two section $\xi_1\in{\h}(\mu_1)$ and $\xi_2\in{\h}(\mu_2)$ which $p_1(x)=\langle L_1(x)\xi_1(d(x)),\xi_1((r(x))\rangle$ and $p_2(x)=\langle L_2(x)\xi_2(d(x)),\xi_2((r(x))\rangle$. Since
$\alpha\phi_{p_1}\leq\phi_p$, $\pi_{L_1}$ are unitary
equivalent to a subrepresentation of $\pi_L$, hence $\mu\sim\mu_1$,
similarly  $\mu\sim\mu_2$. Now if we choice two section $\eta_1,
\eta_2\in{\h}(\mu)$ with $\xi_1\sim\eta_1$ and $\xi_2\sim\eta_2$,
then by Lemma 3.6 and 3.8, $q_1(x)=\langle
L(x)\eta_1(d(x)),\eta_1((r(x))\rangle$ and $q_2(x)=\langle
L(x)\eta_2(d(x)),\eta_2((r(x))\rangle$ which are two elements of
${\mP}^{\mu}_1(G)$  satisfing in $\phi_{p_1}=\phi_{q_1}$ and $\phi_{p_2}=\phi_{q_2}$. Therefore
$\phi_p=\alpha\phi_{q_1}+\beta\phi_{q_2}=\phi_{\alpha q_1+\beta
q_2}$. But $p, \alpha q_1+\beta q_2\in P_1^{\mu}(G)$, so in this case
Lemma 3.6 implies $p=\alpha q_1+\beta q_2$, and the assumption of
$p$ implies $p=q_1=q_2$. Hence $\phi_1=\phi_2=\phi_p$.}
\end{proof}

\begin{thm} Let $\phi\in\Sigma(C^*(G))$  be an extreme point of $\Sigma(C^*(G))$, then there exists $\mu\in Q$ and an extreme point
 $p$ of ${\mP}^{\mu}_1(G)$ such that $\phi=\phi_p$.
\end{thm}
\begin{proof}
{\rm By Lemma 3.5 every element $\phi$ of $\Sigma(C^*(G))$ is of the form $\phi_p$ for some $\mu\in Q$ and some $p\in{\mP}^{\mu}_1(G)$. Assume there are $p_1, p_2\in{\mP}^{\mu}_1(G)$ such that $p=\alpha p_1+\beta p_2$, where $\alpha, \beta\geq 0$ and $\alpha+\beta=1$, then
$\phi_p=\alpha\phi_{p_1}+\beta\phi_{p_2}$. By Remark 3.4 $\phi_{p_1}, \phi_{p_2}\in\Sigma(C^*(G))$, so $\phi_{p}=\phi_{p_1}=\phi_{p_2}$. Note that $p, p_1, p_2\in{\mP}^{\mu}_1(G)$, in this case by Lemma 3.6 $p=p_1=p_2$.}
\end{proof}

\begin{cor} Let $G$ be a locally compact groupoid with Haar system $\lambda$ and let $\mu\in Q$ be an ergodic measure. If
$p\in{\mP}^{\mu}_1(G)$ and $p\in Ext({\mP}^{\mu}_1(G))$, then the representation $L$ which is related to $p$ is an irreducible representation of $G$.
\end{cor}
\begin{proof}
{\rm Suppose $p\in Ext({\mP}^{\mu}_1(G))$, then by theorem 3.9 $\phi_p\in Ext(\Sigma(C^*(G)))$. Hence by [1, 32.7] $\phi_p$ is a pure state on $C^*(G)$ and therefore $\pi_L$ is an irreducible representation of $C^*(G)$ so theorem 2.9 implies $L$ is irreducible}
\end{proof}
\begin{cor} Let $\mu\in Q$ be an ergodic measure, $L$ be an irreducible representation of $G$ on a Hilbert bundle $(G^0,\{H_u\},\mu)$ and $\xi\in{\h}(\mu)$ with $\|\xi\|=1$, then $p(x)=\langle L(x)\xi(d(x)), \xi(r(x))\rangle$ is an extreme point of ${\mP}^{\mu}_1(G)$.
\end{cor}
\begin{proof}
Since $L$ is irreducible, $\pi_L$ is irreducible and by [1, 32.6] $\phi_p$ is a pure state and consequently is an extreme point of $\Sigma(C^*(G))$. By Theorem 3.10 there exists a $\mu'\in Q$ and an extreme point $q$ of ${\mP}^{\mu'}_1(G)$ which $\phi_p=\phi_q$. Therefore $p\sim q$ and it is easy to show that $p\in Ext({\mP}^{\mu}_1(G))$.
\end{proof}


\vskip0.2in
\no {\bf References}
\vskip0.1in

\footnotesize

\REF{[1]}  J. B. Conway, {\it A Course in Functional Analysis},
Springer-Verlag, New York, 1990.

\REF{[2]} J. Dixmier, {\it von Neumann Algebra}, Par II,
North-Holland Publishing Company, Amsterdam, 1981.

\REF{[3]} P. Hahn, {\it Haar measure fore measure groupoids}, Trans. Amer. Math Soc, 242 (1978), 1-33.

\REF{[4]} E. Hewitt and K.A. Ross, Abstruct Harmonic Analysis I, Springer-Verlag, 1963.

\REF{[5]} P. S. Muhly, {\it Coordinates in operator Algebra}, to
appear, CBMS Regional Conference Series in Mathematics, American
Mathematical Society, Providence, 180pp.

\REF{[6]} A. L. T. Paterson, {\it Groupoids, inverse semigroups
and their operator algebras}, Progress in Mathematics, Vol.
$\mathbf{170}$, Birkh$\ddot{a}$user, Boston, 1999.

\REF{[7]} A. L. T. Paterson, {\it The Fourier algebra for locally
compact groupoids, preprint}, 2001.

\REF{[8]} J. Renault. {\it A groupoid approach to $C^*-$
Algebra,} Lecture Note in Mathematics, No.793. Springer-Verlag,
New York-Heidelberg, 1980.

\REF{[9]} J. Renault. {\it The Fourier algebra of a measured groupoid and its multipliers,} J. Functional analysis 145(1997), 455-490.

\REF{[10]} A. Ramsay and M. E. Walter, {\it Fourier-Stieltjes algebra of locally compact groupoids,} J. Functional analysis 148(1997), 314-365.

\end{document}